\documentclass[11pt,reqno]{amsart}
\usepackage{amsmath, amssymb, amsthm}
\usepackage{graphicx}

\numberwithin{equation}{section}

\newtheorem{theorem}{Theorem}[section]

\newtheorem{proposition}[theorem]{Proposition}
\newtheorem{corollary}[theorem]{Corollary}
\newtheorem{lemma}[theorem]{Lemma}

\newtheorem*{general Gromov'}{Corollary \ref{general Gromov}$'$}

\def \proof {\noindent {\bf Proof.}\ \ }

\def \endproof {{\mbox{}\nolinebreak\hfill\rule{2mm}{2mm}\par\medbreak}}

\def \R {\mathbb{R}}
\def \C {\mathbb{C}}

\def \Z {\mathbb{Z}}

\def \F {\mathbb{F}}

\def \a {\alpha}

\def \g {\gamma}
\def \e {\varepsilon}

\def \s {\sigma}

\def \< {\langle}
\def \> {\rangle}

\def \sign {{\rm sign}}
\def \dist {{\rm dist}}

\def \Span {{\rm span}}

\def \range {{\rm range }}

\def \lin {{\rm lin}}
\def \aff {{\rm aff}}

\def \Prob {{\rm Prob}}

\def \supp {{\rm supp}}

\begin{document}
\title [Geometric approach to error correcting codes and signal recovery]
       {Geometric approach to error correcting codes
        and reconstruction of signals}

\author{Mark Rudelson}
\address{Departent of Mathematics, University of Missouri, Columbia, MO 65211, U.S.A.}
\email{rudelson@math.missouri.edu}

\author{Roman Vershynin}
\address{Departent of Mathematics, University of California, Davis, CA 95616, U.S.A.}
\email{vershynin@math.ucdavis.edu}

\thanks{The first author is partially supported by the NSF grant DMS 0245380.
  The second author is partially supported by the NSF grant DMS 0401032
  and by the Miller Scholarship from the University of
  Missouri-Columbia. }


\begin{abstract}
We develop an approach through geometric functional analysis
to error correcting codes and to reconstruction of signals
from few linear measurements. An error correcting code encodes
an $n$-letter word $x$ into an $m$-letter word $y$
in such a way that $x$ can be decoded correctly when any $r$ letters
of $y$ are corrupted. We prove that most linear orthogonal
transformations $Q : \R^n \to \R^m$ form efficient and robust robust
error correcting codes over reals. The decoder (which corrects the corrupted
components of $y$) is the metric projection onto the range of $Q$
in the $\ell_1$ norm. An equivalent problem arises in signal processing:
how to reconstruct a signal that belongs to a small class from few linear measurements?
We prove that for most sets of Gaussian measurements, all signals
of small support can be exactly reconstructed by the $L_1$ norm
minimization. This is a substantial improvement of recent results of Donoho and
of Candes and Tao. An equivalent problem in combinatorial geometry
is the existence of a polytope with fixed number of facets and maximal 
number of lower-dimensional facets.
We prove that most sections of the cube form such polytopes.
\end{abstract}

\maketitle

\section{Error correcting codes and transform coding}

Error correcting codes are used in modern technology to protect
information from errors. Information is formed by finite words
over some alphabet $\F$.
An encoder transforms an $n$-letter word $x$ into an $m$-letter word $y$ with $m > n$.
The decoder must be able to recover $x$ correctly when up to $r$ letters of $y$
are corrupted in any way. Such an encoder-decoder pair is called an
{\em $(n,m,r)$-error correcting code}.

Development of algorithmically efficient error correcing codes 
has been attracting attention of engineers, computer scientists
and applied mathematicians for past five decades.
Known constructions involve deep algebraic and combinatorial methods, 
see \cite{Handbook}, \cite{Sp1}, \cite{Sp2}.
This paper develops a new approach to error correcting codes
from the viewpoint of geometric functional analysis (asymptotic convex geometry).
Our main focus will be on words over the alphabet $\F = \R$ or $\C$. In applications,
these words may be formed of the coefficients of some signal (such as image or audio)
with respect to some basis or overcomplete system (Fourier, wavelet, etc.)
Finite alphabets will be discussed in Section \ref{s:conclusion}.

The simplest and most natural way to encode a vector $x \in \R^n$ into
a vector $y \in \R^m$ is of course a linear transform
\begin{equation}                \label{Q}
y = Qx
\end{equation}
where $Q$ is given by an $m \times n$ matrix. Elementary linear
algebra tells us that if $m \ge n + 2r$ and the range of $Q$ is
generic\footnote{that is, in general position with respect to all
subspaces $\R^I$, $|I| = r$} then $x$ can be recovered from $y$
even if $r$ coordinates of $y$ are corrupted. This gives an
$(n,m,r)$-error correcting code. However, the decoder for this
code has a huge computational complexity, as it involves a search
through all $r$-element subsets of the components of $y$. Then the
problem is:

\medskip

\begin{quote}
{\em How to reconstruct a vector $y$ in an $n$-dimensional subspace $Y$
  of $\R^m$ from a vector $y' \in \R^m$
  that differs from $y$ in at most $r$ coordinates?}
\end{quote}

\medskip

\noindent
What complicates this problem is the arbitrary magnitude of errors in each
corrupted component of $y'$, in contrast to what happens over finite alphabets
such as $\F = \{0,1\}$.

A traditional and simple approach to denoising $y'$, used in applications such
as signal processing, is the mean least square (MLS) minimization. One hopes
that $y$ is well approximated by a solution to the minimization problem
\begin{equation*}
  \min_{u \in Y} \|u - y'\|_2           \tag{MLS}
\end{equation*}
where $\|x\|_2^2 = \sum_i |x_i|^2$.
The solution to (MLS) is simply the orthogonal projection of $y'$ onto $Y$.
This of course can not recover $y$ exactly, and even the approximation is typically
poor since we have no control of the magnitude of the errors in the
corrupted coordinates.
A promising alternative approach is the {\em Basis Pursuit} (BP).
We simply replace the $1$-norm by the $2$-norm and expect $y$ to be the {\em exact}
and unique solution to the minimization problem
\begin{equation*}
  \min_{u \in Y} \|u - y'\|_1           \tag{BP}
\end{equation*}
where $\|x\|_1 = \sum_i |x_i|$.
Thus a solution to (BP) is the metric projection of $y'$ onto $Y$
with respect to the $1$-norm.  (BP) be cast as a Linear Programming problem,
and can be attacked with a variety of methods, such as the classical simplex method
or more recent interior point methods that yield polynomial time algorithms
\cite{CDS}.

\begin{center}
\raisebox{-1 true in}{\includegraphics[height=2in]{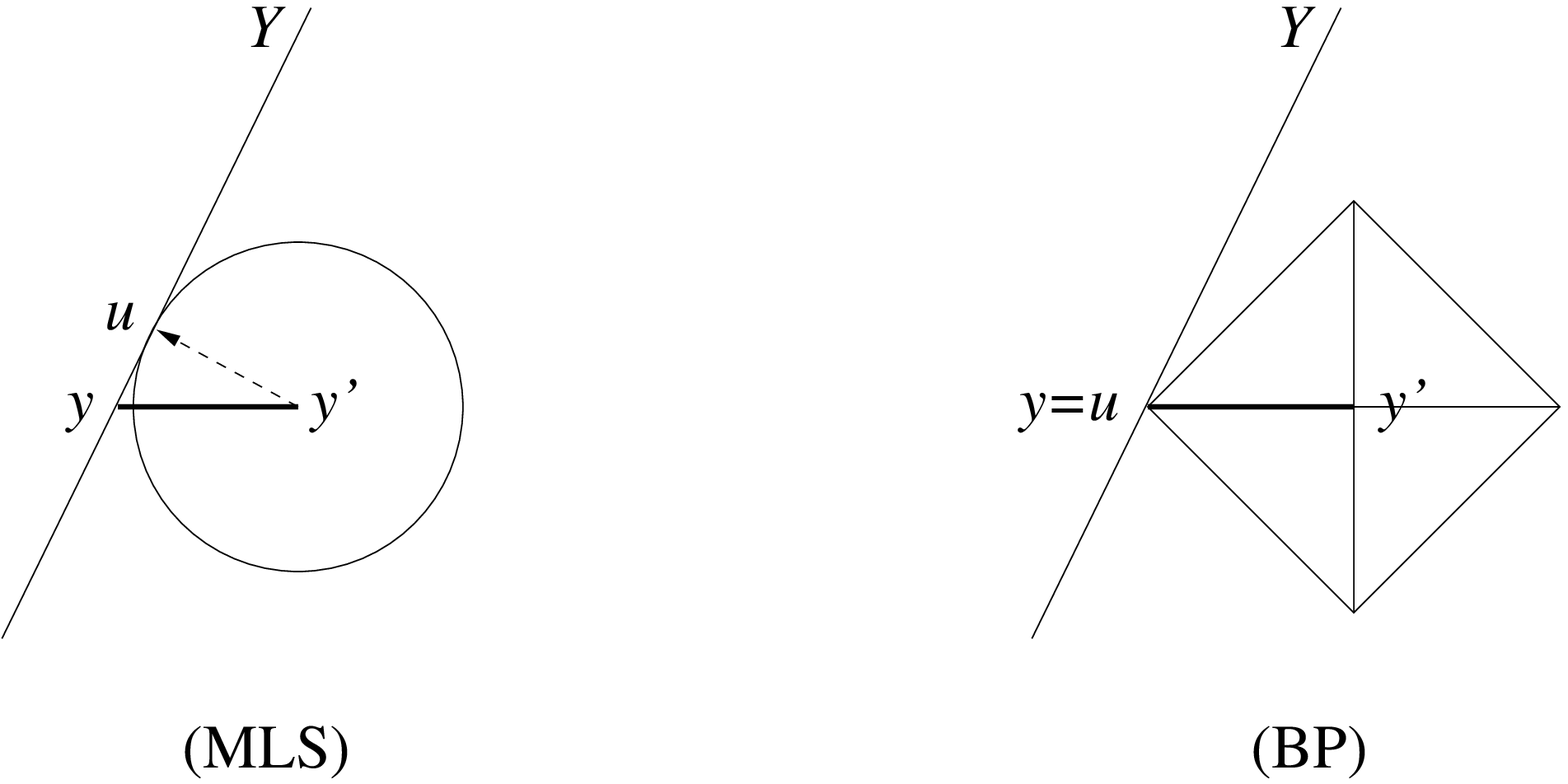}}
\end{center}

The potential of Basis Pursuit for exact reconstruction
is illustrated by the following heuristics, essentially due to \cite{DET}.
The solution $u$ to (MLS) is the contact point where the smallest Euclidean ball
centered at $y'$ meets the subspace $Y$. That contact point is in general
different from $y$. The situation is much better in (BP): typically the solution
coincides with $y$. The solution $u$ to (BP) is the contact point
where the smallest octahedron centered at $y'$ (the ball with respect to the $1$-norm)
meets $Y$. Because the vector $y-y'$ lies in a low-dimensional coordinate subspace,
the octahedron has a wedge at $y$. Thus, many subspaces $Y$ through $y$
will miss the octahedron of radius $y-y'$ (as opposed to the Euclidean ball).
This forces the solution $u$ to (BP), which is the contact point of the octahedron,
to coincide with $y$.

The idea of using the $1$-norm instead of the $2$-norm for better data recovery
has been explored since mid-seventies in various applied areas, in particular
geophysics and statistics (early history can be found in \cite{T 04c}).
With the subsequent development of fast interior point methods in Linear Programming,
(BP) turned into an effectively solvable problem, and was put forward
more recently by Donoho and his collaborators, triggering
massive experimental and theoretical work \cite{CDS, DH, EB, FN, DE, GN,
T 04a, T 04b, T 04c, DET, D 04a, D 04b, DT 04a, DT 04b, CRT, CR, CT}.

\medskip

The main result of this paper validates the Basis Pursuit method
for most subspaces $Y$ under an asymptotically sharp condition on $m,n,r$.
We thus prove that {\em the Basis Pursuit yields exact reconstruction for most subspaces $Y$}
in the Grassmanian.
The randomness is with respect to the normalized Haar
measure on the Grassmanian $G_{m,n}$ of $n$-dimensional subspaces of $\R^m$.
Positive absolute constants will be denoted throughout the paper
by $C, c, C_1, \ldots$.

\begin{theorem}                                 \label{ecc}
  Let $m$, $n$ and $r < cm$ be positive integers such that
  \begin{equation}              \label{mnr'}
    m = n+ R, \ \ \ \text{where $R \ge C r \log(m/r)$}.
  \end{equation}
  Then a random $n$-dimensional subspace $Y$ in $\R^m$ satisfies
  the following with probability at least $1 - e^{-c R}$.
  Let $y \in Y$ be an unknown vector, and we are given a vector $y'$ in $\R^m$
  that differs from $y$ on at most $r$ coordinates.
  Then $y$ can be exactly reconstructed from $y'$ as the solution
  to the minimization problem (BP).
\end{theorem}

In an equivalent form, this theorem is a substantial improvement of
recent results of Donoho \cite{D 04a} and of Candes and Tao \cite{CT},
see Theorem~\ref{reconstruction} below.

\subsection{Error correcting codes.}                \label{ss:ecc}
Theorem \ref{ecc} implies a natural $(n,m,r)$-error correcting code over $\R$.
The encoder \eqref{Q} is given by an $m \times n$
random orthogonal matrix\footnote{one can view it as the first $n$ rows of
a random matrix from $O(m)$ equipped with the normalized Haar measure.} $Q$.
Its range $Y$ is a random $n$-dimensional subspace in $\R^m$.
The decoder takes a corrupted vector $y'$, solves (BP) and outputs
$Q^T u = Q^{-1} u$. Theorem \ref{ecc} states that under the assumption \eqref{mnr'},
this encoder-decored pair is an $(n,m,r)$-error correcting code with
exponentially good probability $\ge 1 - e^{-c R}$.

\subsection{Sharpness.}
The sufficient condition \eqref{mnr'} is sharp up to an absolute
constant $C$ (see Section \ref{s:conclusion}) and is only slightly
stronger than the necessary condition $m \ge n + 2r$. The ratio
$\e = r/m$ in \eqref{mnr'} is the number of errors per letter in
the noisy communication channel that maps $y$ to $y'$. Thus $\e$
should be considered as a quality of the channel, which is
independent of the message. Thus \eqref{mnr'} is equivalent to
$$
m \ge \Bigl(1 + C \e \log \frac{1}{\e} \Bigr) n.
$$

\subsection{Robustness.}
An natural feature of our error correction code is its {\em robustness}.
Simple linear algebra yields that
the solution to (BP) is stable with respect to the $1$-norm -- in the same way
as the solution to (MLS) is stable with respect to the $2$-norm, see \cite{CT}.
Such robustness allows in particular quantization of the messages.
This immediately yields error correcting codes for finite alphabets, see
Section \ref{s:conclusion}.

\subsection {Transform coding.}
In the signal processing, the linear codes \eqref{Q} are known
as {\em transform codes}. The general paradigm about transform codes is
that the redundancies in the coefficients of $y$
that come from the excess of the dimension $m > n$ should guarantee
a stability of the signal with respect to noise, quantization, erasures,
etc. This is confirmed by an extensive experimental and some theoretical
work, see e.g. \cite{Da,G1,G2,GVT,GKK,KDG,BO,CK}
and the bibliography contained therein.
Theorem \ref{ecc} states that {\em most orthogonal transform codes
are good error-correcting codes}.

\subsection* {Acknowledgement.} This work has started when the second
author was visiting University of Missouri-Columbia as a Miller
Visiting Scholar. He is grateful to UMC for the hospitality.

\section{Reconstruction of signals from linear measurements.}

The heuristic idea that guides the Statistical Learning Theory is that
{\em a function $f$ from a small class should be determined by few linear measurements}.
Linear measurements are generally given by some linear functionals $X_k$
in the dual space, which are fixed (in particular are independent of $f$).
Most common measurements are point evaluation functionals; the
problem there is to interpolate $f$ between known values while keeping $f$
in the known (small) class.
When the evaluation points are chosen at random, this becomes the `proper learning'
problem of the Statistical Learning Theory (see \cite{M}).

We shall however be interested in general linear measurements.
The proposal to learn $f$ from general linear measurements ({\em `sensing'})
has been originated recently from a criticism of the current methodology
of signal compression. Most of real life signals, such as images and sounds,
seem to belong to small classes. This is because they carry much of unwanted information
that can be discarded with almost no perceptual loss, which makes such signals
easily compressible. Donoho \cite{D 04c} then questions the conventional scheme of
signal processing, where the whole signal must be first acquired (together
with lots of unwanted information) and only then be compressed
(throwing away the unwanted part).
Instead, can one {\em directly acquire} (`sense') the essential part of the signal,
via few linear measurements? Similar issues are raised in \cite{CT}.
We shall operate under the assumption that some technology
allows us to take linear measurements in certain fixed `directions' $X_k$.

We will assume that our signal $f$ is discrete, so we view it as a vector in $\R^m$.
Suppose we can take linear measurements $\< f, X_k \> $ with some fixed
vectors $X_1, X_2, \ldots, X_{R}$ in $\R^m$.
Assuming that $f$ belongs to a small class,
how many measurements $R$ are needed to reconstruct $f$?
And even when we prove that $R$ measurements do determine $f$
(uniquely or approximately), the algorithmic issue remains unsettled:
how can one reconstruct $f$ from these measurements?

The previoous section suggests to reconstruct $f$
as a solution to the Basis Pursuit minimization problem
\begin{equation*}
  \min \|g\|_1
  \ \ \text{subsect to} \ \
  \< g, X_k \> = \< f, X_k \> , \ \ k = 1, \ldots, R.     \tag{BP$'$}
\end{equation*}
For the Basis Pursuit to work, the vectors $X_k$ must be in a good position
with respect to all coordinate subspaces $\R^I$, $|I| \le r$.
A typical choice for such vectors would be the independent standard Gaussian
vectors\footnote{All the components of $X_k$ are independent
standard Gaussian random variables.} $X_k$.

\subsection{Functions with small support}
In the class of functions with small support, one can hope for exact reconstruction.
Candes and Tao \cite{CT} have indeed proved that every {\em fixed} function $f$ with
support $|\supp f| \le r$ can indeed be recovered by (BP$'$), correctly
with the polynomial probability $1 - m^{-\text{const}}$, from the
$R = C r \log m$ Gaussian measurements.
However, the polynomial probability is clearly not sufficient
to deduce that there is {\em one} set vectors $X_k$ that can be used to
reconstruct all functions $f$ of small support.

The following equivalent form of Theorem \ref{ecc} does
yield a uniform exact reconstruction.
It provides us with {\em one set} of linear measurements from from which we
can effectively reconstruct {\em every} signal of small support.

\begin{theorem} [Uniform Exact Reconstruction]                   \label{reconstruction}
  Let $m$, $r < cm$ and $R$ be positive integers satisfying
  $R \ge C r \log(m/r)$.
  The independent standard Gaussian vectors $X_k$ in $\R^m$
  satisfy the following with probability at least $1 - e^{-c R}$.
  Let $f \in \R^m$ be an unknown function of small support, $|\supp f| \le r$,
  and we are given $R$ measurements $\< f, X_k\> $.
  Then $f$ can be exactly reconstructed from these measurements
  as a solution to the Basis Pursuit problem (BP$'$).
\end{theorem}

This theorem gives uniformity in Candes-Tao result \cite{CT}, improves the polynomial
probability to an exponential probability, and improves upon the number $R$
of measurements (which was $R \ge C r \log m$ in \cite{CT}).  
Donoho \cite{D 04c} proved a weaker form of Theorem \ref{reconstruction}
with $R/r$ bounded below by some function of $m/r$.

\medskip

\proof
Write $g = f - u$ for some $u \in \R^m$. Then (BP$'$) reads as
\begin{equation}                    \label{BP uf}
\min \|u - f\|_1
\ \ \text{subsect to} \ \
\< u, X_k \> = 0, \ \ k = 1, \ldots, R.
\end{equation}
The constraints here define a random $(n = m - R)$-dimensional subspace
$Y$ of $\R^m$. Now apply Theorem \ref{ecc} with $y = 0$ and $y' = f$. It states
that the unique solution to \eqref{BP uf} is $u = 0$. Therefore, the
unique solution to (BP$'$) is $f$.
\endproof

\subsection{Compressible functions}
In a larger class of compressible functions \cite{D 04c}, we can only hope for
an approximate reconstruction. This is a class of functions $f$ that are
well compressible by a known orthogonal transform, such as Fourier or wavelet.
This means that the coefficients of $f$ with respect to a certain known
orthogonal basis have a power decay. By applying an appropriate rotation,
we can assume that this basis is the canonical basis of $\R^m$, thus
$f$ satisfies
\begin{equation}                        \label{compressible}
  f^*(s) \le s^{-1/p}, \ \ \ s = 1, \ldots, m
\end{equation}
where $f^*$ denotes a nonincreasing rearrangement of $f$.
Many natural signals are compressible for some $0 < p < 1$,
such as smooth signals and signals with bounded variations (see \cite{CT}),
in particular most photographic images.
Theorem \ref{reconstruction} implies, by the argument of \cite{CT},
that functions compressible in some basis can be approximately
reconstructed from few fixed linear measurements:

\begin{corollary}[Uniform Approximate Reconstruction]
  Let $m$ and $r$ be positive integers.
  The independent standard Gaussian vectors $X_k$ in $\R^m$
  satisfy the following with probability at least $1 - e^{-c R}$.
  Assume that an unknown function $f \in \R^m$ satisfies either
  \eqref{compressible} for some $0 < p < 1$ or $\|f\|_1 \le 1$ for $p=1$.
  Suppose that we are given $R$ measurements $\< f, X_k\> $.
  Then $f$ can be approximately reconstructed from these measurements:
  a unique solution $g$ to the Basis Pursuit problem (BP$'$) satisfies
  $$
  \|f - g\|_2
    \le C_p \Bigl( \frac{\log(m/R)}{R} \Bigr)^{\frac{1}{p} - \frac{1}{2}}
  $$
  where $C_p$ depends on $p$ only.
\end{corollary}

This theorem also gives uniformity in another Candes-Tao result from \cite{CT}
(see also \cite{D 04b}); it improves the polynomial probability to an
exponential probability, and also improves upon the approximation error.

\section{Counting low-dimensional facets of polytopes.}

Theorem \ref{ecc} turns out to be equivaent to a problem of counting 
lower-dimensional facets of polytopes. Let $B_1^m$ denote the unit ball 
with respect to the $1$-norm; it is sometimes called the unit octahedron.
The polar body is the unit cube $B_\infty^m = [-1,1]^m$. 
The conclusion of Theorem \ref{ecc} is then equivalent to the
following statement: the affine subspace $z + Y$ is tangent to the unit
octahedron at point $z$, where $z = y' - y$. This should happen
for all $z$ from the coordinate subspaces $\R^I$ with $|I| = r$.
By the duality, this means that the subspace $Y^\perp$ intersects all
$(m-r)$-dimensional facets of the unit cube.  The section of the cube by
the subspace $Y^\perp$ forms an origin-symmetric polytope of dimension $R$
and with $2m$ facets.

Our problem can thus be stated as a problem of counting lower-dimensional facets
of polytopes.
\begin{quote}
  {\em Consider an $R$-dimensional origin symmetric polytope
  with $2m$ facets. How many $(R-r)$-dimensional facets can it have?}
\end{quote}
Clearly\footnote{Any such facet is the intersection of some $r$ facets
of the polytope of full dimension $R-1$; there are $m$ facets to choose from,
each coming with its opposite by the symmetry.}, no more than
$2^r \binom{m}{r}$. Does there exist a polytope with that many facets?
Our ability to construct such a polytope
is equivalent to the existence of the efficient error
correcting code. Indeed, looking at the canonical realization of such a
polytope as a section of the unit cube by a subspace $Y^\perp$,
we see that $Y^\perp$ intersects all the $(m-r)$-dimensional facets
of the cube. Thus $Y$ satisfies the conclusion of Theorem~\ref{ecc}.
We can thus state Theorem \ref{ecc} in the following form:

\begin{theorem}
  There exists an $R$-dimensional symmetric polytope with $m$ facets
  and with the maximal number of $(R-r)$-dimensional facets
  (which is $2^r \binom{m}{r}$), provided $R \ge C r \log(m/r)$.
  A random section of the cube forms such a polytope with probability
  $1 - e^{-cR}$.
\end{theorem}

So, how can we prove that a random subspace $Y^\perp$ indeed intersects all the
$(m-r)$-dimensional facets of the cube? It is enough to show that
$Y^\perp$ intersects one such fixed facet with exponential probability
(bigger than $1 - 2^{-r} \binom{m}{r}^{-1}$).
The main difficulty here is that the concentration of measure technique
can not be readily applied. This is because the $\infty$-norm defined
by the unit cube (more precisely, by its facet) has a bad Lipschitz constant.
To improve the Lipschitzness, we first project the facet onto a random
subspace (within its affine span); the random subspace parallel to which we
project is taken from the random directions that form $Y^\perp$.
This creates a big Euclidean ball inside the projected facet;
here we shall use the full strength of the estimate
of Garnaev and Gluskin \cite{GG} on Euclidean projections of a cube.
The existence of the Euclidean ball inside a body creates the needed
Lipschitzness, so we can now use the concentration of measure tecnique.

\medskip

The rest of the paper is organized as follows.
In Section \ref{s:proof} we prove Theorem \ref{ecc}.
In Section \ref{s:conclusion} we discuss some optimality and
robustness of the Basis Pursuit with applications to error correcting
codes over finite alphabets.

\section{Proof}                         \label{s:proof}

We shall use the following standard notations throughout the proof.
The $p$-norm ($1 \le p < \infty$) on $\R^m$ is defined by
$\|x\|_p^p = \sum_i |x_i|^p$, and for $p = \infty$ it is
$\|x\|_\infty = \max_i |x_i|$. The unit ball with respect to the
$p$-norm on $\R^n$ is denoted by $B_p^m$. When the $p$-norm is considered
on a coordinate subspace $\R^I$, $I \subset \{1,\ldots,m\}$,
the corresponding unit ball is denoted by $B_p^I$.

The unit Euclidean sphere in a subspace $E$ is denoted by $S(E)$.
The normalized rotational invariant Lebesgue measure on $S(E)$ is denoted
by $\sigma_E$.
The orthogonal projection in onto a subspace $E$ is denoted by $P_E$.
The standard Gaussian measure on $E$ (with the identity covariance matrix)
is denoted by $\gamma_H$. When $E = \R^d$, we write $\sigma_{d-1}$ for
$\sigma_E$ and $\gamma_d$ for $\gamma_E$.

\subsection{Duality}
We begin the proof of Theorem \ref{ecc} with a typical duality argument,
leading to the same reformulation of the problem as in \cite{CT}.
We claim that the conclusion of Theorem \ref{ecc} follows from
(and is actually equivalent to) the following separation condition:
\begin{equation}                        \label{separation}
  (z + Y) \cap \;\text{interior}\, (B_1^m) = \emptyset
  \ \ \ \text{for all} \ \ z \in \bigcup_{|I| = r} B_1^I.
\end{equation}
Indeed, suppose \eqref{separation} holds. We apply it for
$$
z := \frac{y-y'}{\|y-y'\|_1}
$$
noting that $z \in \bigcup_{|I| = r} B_1^I$ holds, because $y$ and $y'$
differ in at most $r$ coordinates.
By \eqref{separation},
$$
(z + v) \cap \;\text{interior}\, (B_1^m) = \emptyset
\ \ \ \text{for all $v \in Y$}
$$
which implies
$$
\|z + v\|_1 \ge 1
\ \ \ \text{for all $v \in Y$}.
$$
Let $u \in Y$ be arbitrary. Using the inequality above for
$v := \frac{u-y}{\|u-y\|_1}$, we conclude that
$$
\|u-y\|_1 \ge \|y-y'\|_1
\ \ \ \text{for all $u \in Y$}.
$$
This proves that $y$ is indeed a solution to (BP).
The solution to (BP) is unique with probability $1$ in the Grassmanian.
This follows from a direct dimension argument, see e.g. \cite{CT}.

By Hahn-Banach theorem, the separation condition \ref{separation}
is equivalent to the following:
for every $z \in \bigcup_{|I| = r} \;\text{boundary}\, B_1^I$
there exists $w = w(z) \in Y^\perp$ such that
$$
\< w,z \> = \sup_{x \in B_1^m} \< w,x \> = \|w\|_\infty.
$$
This holds if and only if the components of $w$ satisfy
\begin{equation}                        \label{w}
  \begin{cases}
    w_j = \sign(z_j) \ \ \text{for $j \in I$}, \\
    |w_j| \le 1 \ \ \text{for $j \in I^c$}.
  \end{cases}
\end{equation}
The set of vectors $w$ in $\R^m$ that satisfy \eqref{w} form a
$(m-r)$-dimensional facet of the unit cube $B_\infty^m$.
Then with $E := Y^\perp$ we can say that the conclusion
of Theorem \ref{ecc} is equivalent to the following:

\medskip

\begin{quote}
  {\em A random $R$-dimensional subspace $E$ in $\R^m$ intersects
   all the $(m-r)$-dimensional facets of the unit cube
   with probability at least $1 - e^{-cR}$.}
\end{quote}

\medskip

It will be enough to show that $E$ intersects {\em one fixed}
facet with the probability $1 - e^{-cR}$. Indeed, since the total
number of the facets is $N = 2^r \binom{m}{r}$, the probability
that $E$ misses some facet would be at most $N e^{-cR} \le e^{-c_1 R}$
with an appropriate choice of the absolute constant in \eqref{mnr'}.

\subsection{Realizing a random subspace}
We are to show that a random $R$-dimensional subspace $E$ intersects one fixed
$(m-r)$-dimensional facet of the unit cube $B_\infty^m$ with high probability.
Without loss of generality, we can assume that our facet is
$$
F = \{ (w_1, \ldots, w_{m-r}, 1, \ldots, 1), \ \ \text{all $|w_j| \le 1$} \},
$$
whose center is
$$
\theta = (\underbrace{0,\ldots,0}_{m-r}, 1,\ldots,1).
$$
The probability we are interested in is
$$
P := \Prob\{ E \cap F \ne \emptyset\}.
$$
We shall restrict our attention to the linear span of $F$,
$$
\lin(F) = \{ (w_1, \ldots, w_{m-r}, t, \ldots, t),
     \ \ \text{all $w_j \in \R$, $t \in \R$} \},
$$
and even to its the affine span of $F$,
$$
\aff(F) = \{ (w_1, \ldots, w_{m-r}, 1, \ldots, 1),
     \ \ \text{all $w_j \in \R$} \}.
$$
Only the random affine subspace $E \cap \aff(F)$ matters for us, because
$$
P =  \Prob\Bigl\{ (E \cap \aff(F)) \cap F \ne \emptyset \Bigr\}.
$$
The dimension of that affine subspace is almost surely
$$
l := \dim (E \cap \aff(F)) = R-r.
$$

We can realize the random affine subspace $E \cap \aff(F)$
(or rather a random subspace with the same law) by the following
algorithm:

\begin{enumerate}

  \item Select a random variable $D$ with the same law as
    $\dist(\theta, E \cap \aff(F))$.

  \item Select a random subspace $L_0$ in the Grassmanian $G_{m-r,l}$.
    It will realize the ``direction'' of $E \cap \aff(F)$ in $\aff(F)$.

  \item Select a random point $z$ on the Euclidean sphere $D \cdot S(L_0^\perp)$
    of radius $D$, according to the uniform distribution on the sphere.
    Here $L_0^\perp$ is the orthogonal complement of $L_0$ in $\R^{m-r}$.
    The vector $z$ will realize the distance from the affine subspace
    $E \cap \aff(F)$ to the center $\theta$ of $F$.

  \item Set $L = \theta + z + L_0$. Thus the random affine subspace $L$
    has the same law as $E \cap \aff(F)$.

\end{enumerate}

\begin{center}
\raisebox{-1 true in}{\includegraphics[height=2in]{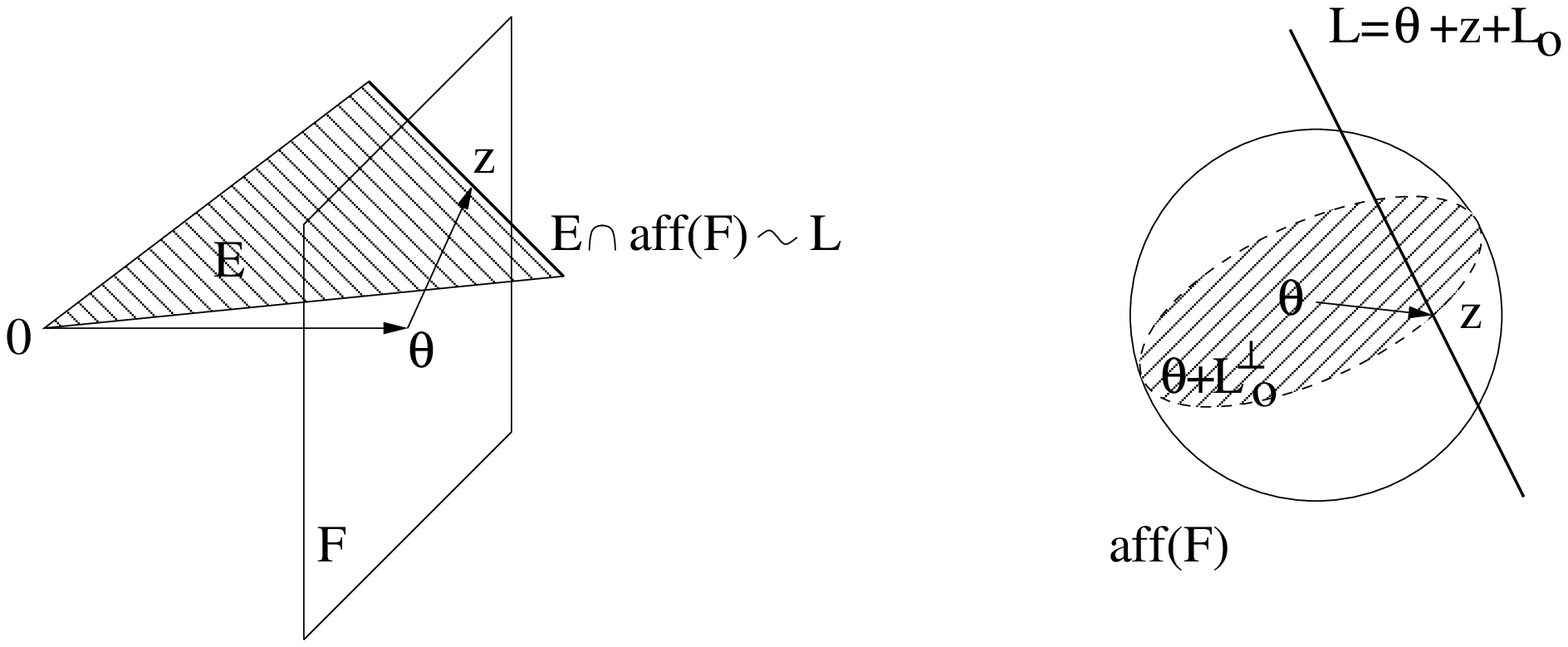}}
\end{center}

\noindent Hence
$$
P = \Prob \{ L \cap F \ne \emptyset \}
  = \Prob \{ (z + L_0) \cap B_\infty^{m-r} \ne \emptyset \}
  = \Prob \{ z \in P_{L_0^\perp} B_\infty^{m-r} \}.
$$
$H := L_0^\perp$ is a random subspace in $G_{m-r,m-r-l} = G_{m-r,m-R}$.
By the rotational invariance of $z \in D \cdot S(H)$,
\begin{equation}                    \label{P=integral}
P = \int_{\R^+} \int_{G_{m-r,m-R}} \sigma_H (D^{-1} P_H B_\infty^{m-r})
      \; d\nu(H) \; d\mu(D)
\end{equation}
where $\nu$ is the normalized Haar measure on $G_{m-r,m-R}$
and $\mu$ is the law of $D$.
We shall bound $P$ in two steps:

\begin{enumerate}

\item Prove that the distance $D$ is small with high probability;

\item Prove that a suitable multiple of the random projection
  $P_H B_\infty^{m-r}$ has an almost full Gaussian
  (thus also spherical) measure.

\end{enumerate}

\subsection{The distance $D$ from the center of the facet to a random subspace}
We shall first relate $D$, the distance to the affine subspace $E \cap \aff(F)$,
to the distance to the linear subspace $E \cap \lin(F)$.
Equivalently, we compute the length of the projection onto $E \cap \lin(F)$.

\begin{lemma}                       \label{linear vs affine}
$$
\|P_{E \cap \lin(F)} \theta \|_2 = \sqrt{\frac{r}{r+D^2}} \;
\|\theta\|_2.
$$
\end{lemma}

\proof
Let $f$ be the multiple of the vector $P_{E \cap \lin(F)} \theta$ such that
$f-\theta$ is orthogonal to $\theta$. Such a multiple exists and is unique,
as this is a two-dimensional problem.

\begin{center}
\raisebox{-1 true in}{\includegraphics[height=1.25in]{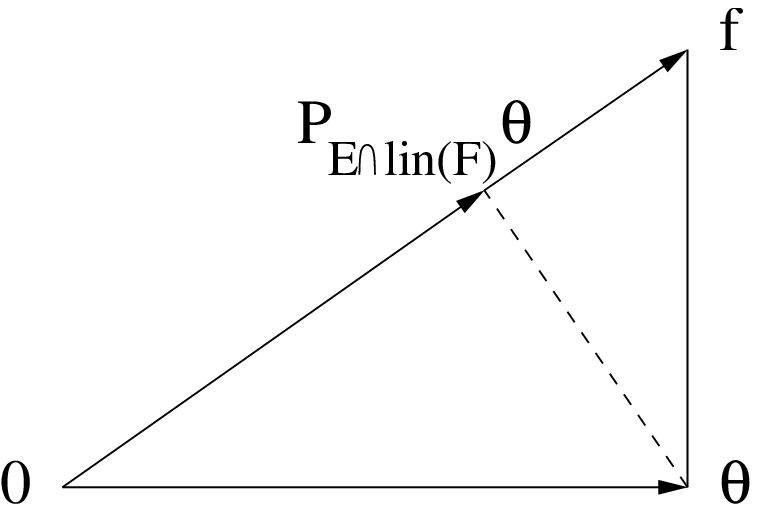}}
\end{center}

Then $f \in E \cap \aff(F)$. Notice that $D= \|f-\theta\|_2$.  By
the similarity of the triangles with the vertices $(0, \theta,
P_{E \cap \lin(F)} \theta)$ and $(0, f, \theta)$, we conclude that
$$
\|P_{E \cap \lin(F)} \theta \|_2 = \frac{r}{\sqrt{r+D^2}} =
\sqrt{\frac{r}{r+D^2}} \; \|\theta\|_2
$$
because $\|\theta\|_2 = \sqrt{r}$.
This completes the proof.
\endproof

\medskip

The length of the projection of a fixed vector onto a random subspace in
Lemma~\ref{linear vs affine} is well known. The asymptotically sharp
estimate was computed by S.~Artstein \cite{A}, but we will be satisfied
with a much weaker elementary estimate, see e.g. \cite{Ma} 15.2.2.

\begin{lemma} \label{l: random projection}
  Let $\theta \in \R^{d-1}$ and let $G$ be a random subspace in $G_{d,k}$.
  Then
  $$
  \Prob \Bigl\{ c \sqrt{\frac{k}{d}} \; \|\theta\|_2
                \le \|P_G \theta\|_2
                \le C \sqrt{\frac{k}{d}} \; \|\theta\|_2
        \Bigr\}
  \ge  1 - 2 e^{-ck}.
  $$
\end{lemma}

We apply this lemma for $G = E \cap \lin(F)$, which is a random subspace
in the Grassmanian of $(l+1)$-dimensional subspaces of $\lin(F)$.
Since $\dim \lin(F) = m-r+1$, we have
$$
\Prob \Bigl\{ \|P_{E \cap \lin(F)} \theta\|_2
        \ge c \sqrt{\frac{l+1}{m-r+1}} \; \|\theta\|_2
      \Bigr\}
\ge  1 - 2 e^{-cl}.
$$
Together with Lemma \ref{linear vs affine} this gives
\begin{equation}                    \label{D small}
\Prob \Bigl\{ D \le c \sqrt{m-r} \sqrt{\frac{r}{l}} \Bigr\}
\ge 1 - 2e^{-cl}.
\end{equation}
Note that $\sqrt{m-r}$ is the radius of the Euclidean ball circumscribed
on the facet $F$. The statement $D \le \sqrt{m-r}$ would only tell us
that the random subspace $E$ intersects the circumscribed ball, not yet the
facet itself. The ratio $r/l$ in \eqref{D small} will be chosen logarithmically
small, which will force $E$ intersect also the facet $F$.

\subsection{Gaussian measure of random projections of the cube}
By \eqref{P=integral} and \eqref{D small},
$$
P \ge \int_{G_{m-r,m-R}}
      \sigma_H \Bigl( \frac{c}{\sqrt{m-r}} \sqrt{\frac{l}{r}} \,
                       P_H B_\infty^{m-r} \Bigr)
      \; d\nu(H) -2 e^{-cl}.
$$
We can replace the spherical measure $\sigma_H$ by the
Gaussian measure $\g_H$ via a simple lemma:

\begin{lemma}                       \label{spherical vs Gaussian}
  Let $K$ be a star-shaped set in $\R^d$. Then
  $$
  \g_d(c \sqrt{d} \cdot K) - e^{-d}
  \le \sigma_{d-1}(K)
  \le \g_d(C \sqrt{d} \cdot K)\cdot (1+ e^{-d}).
  $$
\end{lemma}

\proof Passing to polar coordinates, by the rotational invariance
of the Gaussian measure we see that there exists a probability
measure $\mu$ on $\R^+$ so that the Gaussian measure of every set
$A$ can be computed as $\int_{\R^+} \s^t(A) \; d\mu(t)$, where
$\s^t$ denotes the normalized Lebesgue measure on the Euclidean
sphere of radius $t$ in $\R^d$. Since $K$ is star-shaped,
$\s^t(K)$ is a non-increasing function of $t$. Hence
\begin{align*}
  \gamma_d(K)
  & \ge \int_0^{C \sqrt{d}} \s^t(K) \, d\mu(t)
    \ge \s^{C \sqrt{d}}(K) \cdot
    \gamma_d( C \sqrt{d} B_2^d)
\intertext{and}
  \gamma_d(K)
  & \le \int_0^{c \sqrt{d}} d\mu(t)
   + \s^{c \sqrt{d}}(K) \int_{c \sqrt{d}}^\infty d\mu(t)
  \le \gamma_d(c \sqrt{d} \cdot B_2^d) + \s^{c \sqrt{d}}(K).
\end{align*}
The classical large deviation inequalities imply $\gamma_d(c
\sqrt{d} \cdot B_2^d) \le e^{-d}$ and $\gamma_d( C \sqrt{d} B_2^d)
\ge 1- e^{-d}/2$. Using the above argument for $c \sqrt{d} \cdot
K$, we conclude that $\g_d(c \sqrt{d} \cdot K) \le e^{-d} +
\sigma_{d-1}(K)$ and $\g_d(C \sqrt{d} \cdot K) \ge \sigma_{d-1}(K)
\cdot (1-e^{-d}/2)$.
\endproof

\medskip

Using Lemma \ref{spherical vs Gaussian}
in the space $H$ of dimension $d = m-R$, we obtain
$$
P \ge \int_{G_{m-r,m-R}}
      \gamma_H \Bigl( c \sqrt{\frac{m-R}{m-r}} \sqrt{\frac{l}{r}} \,
                       P_H B_\infty^{m-r} \Bigr)
      \; d\nu(H) -2 e^{-cl} -e^{m-R}.
$$
By choosing the absolute constant $c$ in the assumption $r < cm$
appropriately small, we can assume that $2r < R < m/2$.
Thus
\begin{equation}                        \label{P}
P \ge \int_{G_{m-r,m-R}}
      \gamma_H \Bigl( c \sqrt{\frac{R}{r}} \,
                       P_H B_\infty^{m-r} \Bigr)
      \; d\nu(H) -2 e^{-cR}.
\end{equation}
We now compute the Gaussian measure of random projections of the cube.

\begin{proposition}                     \label{proj of cube}
  Let $H$ be a random subspace in $G_{n,n-k}$, $k < n/2$.
  Then the inequality
  $$
  \gamma_H \Bigl( C \sqrt{\log \frac{n}{k}} \,
                       P_H B_\infty^n \Bigr)
  \ge 1 - e^{-ck}
  $$
  holds with probability at least $1 - e^{-ck}$ in the Grassmanian.
\end{proposition}

The proof of this estimate will follow from the concentration of Gaussian measure,
combined with the existence of a big Euclidean ball inside a random projection
of the cube.

\begin{lemma}[Concentration of Gaussian measure]        \label{concentration}
  Let $A$ be a measurable set in $\R^n$. Then for $\e > 0$,
  $$
  \gamma_n(A) \ge e^{-\e^2 n}
  \ \ \ \text{implies} \ \ \
  \gamma_n(A + C \e \sqrt{n} B_2^n ) \ge 1 - e^{-\e^2 n}.
  $$
\end{lemma}

With the stronger assumption $\gamma(A) \ge 1/2$, this lemma is the classical
concentration inequality, see \cite{L} 1.1. The fact that the concentration
holds also for exponentially small sets follows formally by a simple extension
argument that was first noticed by D.~Amir and V.~Milman in \cite{AM},
see \cite{L} Lemma 1.1.

The optimal result on random projections of the cube
is due to Garnaev and Gluskin \cite{GG}.

\begin{theorem}[Euclidean projections of the cube \cite{GG}]            \label{GG lemma}
  Let $H$ be a random subspace in $G_{n,n-k}$, where $k = \a n < n/2$.
  Then with probability at least $1 - e^{-ck}$ in the Grassmanian, we have
  $$
  c(\a) \, P_H(\sqrt{n} B_2^n)
  \subseteq P_H(B_\infty^n) \subseteq
  P_H(\sqrt{n} B_2^n)
  $$
  where
  $$
  c(\a) = c \sqrt{\frac{\a}{\log(1/\a)}}.
  $$
\end{theorem}

\medskip

\noindent {\bf Proof of Proposition \ref{proj of cube}. }
Let $g_1, g_2, \ldots$ be independent standard Gaussian random variables.
Then for a suitable positive absolute constant $c$ and for every $0 < \e < 1/2$,
$$
\gamma_n \Bigl( C \sqrt{\log \frac{1}{\e}} \, B_\infty^n \Bigr)
= \Prob \Bigl\{ \max_{1 \le j \le n} |g_i| \le C \sqrt{\log \frac{1}{\e}} \Bigr\}
\ge (1 - \e^2/10)^n \ge e^{-\e^2 n}.
$$
Since for every measurable set $A$ and every subspace $H$ one has
$\gamma_H(P_H A) \ge \gamma(A)$, we conclude that
$$
\gamma_H \Bigl( C \sqrt{\log \frac{1}{\e}} \, P_H B_\infty^n \Bigr)
\ge e^{-\e^2 n}
\ \ \ \text{for $0 < \e < 1/2$.}
$$
Then by Lemma \ref{concentration},
\begin{equation}                            \label{cube+ball}
\gamma_H \Bigl( C \sqrt{\log \frac{1}{\e}} \, P_H B_\infty^n
  + C \e \sqrt{n} \, P_H B_2^n \Bigr)
\ge 1 - e^{-\e^2 n}
\ \ \ \text{for $0 < \e < 1/2$.}
\end{equation}
Theorem \ref{GG lemma} tells us that for a random subspace $H$,
if $\e = c \sqrt{\a} = c \sqrt{k/n}$,
then Euclidean ball is absorbed by the projection of the cube
in \eqref{cube+ball}:
$$
\e \sqrt{n} \, P_H B_2^n \subset C \sqrt{\log \frac{1}{\e}} \, P_H B_\infty^n.
$$
Hence for a random subspace $H$ and for $\e$ as above we have
$$
\gamma_H \Bigl( C \sqrt{\log \frac{1}{\e}} \, P_H B_\infty^n \Bigr)
\ge 1 - e^{-\e^2 n},
$$
which completes the proof.
\endproof

\medskip

Coming back to \eqref{P}, we shall use Lemma \ref{proj of cube}
for a random subspace $H$ in the Grassmanian $G_{m-r,m-R}$.
We conclude that if
\begin{equation}                            \label{Rr}
c \sqrt{\frac{R}{r}} \ge C \sqrt{\log \frac{m-r}{R-r}},
\end{equation}
then with probability at least $1 - e^{-cR}$ in the Grassmanian,
$$
\gamma_H \Bigl( c \sqrt{\frac{R}{r}} \, P_H B_\infty^{m-r} \Bigr)
\ge 1 - e^{-cR}.
$$
Since $\frac{m-r}{R-r} \le \frac{m}{r}$, the choice of $R$ in \eqref{mnr'}
satisfies condition \eqref{Rr}. Thus \eqref{P} implies
$$
P \ge 1 - 3 e^{-cR}.
$$
This completes the proof.
\endproof

\section{Optimality, robustness, finite alphabets}                    \label{s:conclusion}

\subsection{Optimality}
The logarithmic term in Theorems \ref{ecc} and
\ref{reconstruction} is necessary, at least in the case of small
$r$.  Indeed, combining  formula \eqref{P=integral} and Lemmas
\ref{linear vs affine}, \ref{l: random projection}, \ref{spherical
vs Gaussian}, we obtain
\begin{equation}  \label{upper P}
  P \le \int_{G_{m-r,m-R}}
      \gamma_H \Bigl( c \sqrt{\frac{R}{r}} \,
                       P_H B_\infty^{m-r} \Bigr)
      \; d\nu(H) + 2 e^{-cR}.
\end{equation}
To estimate the Gaussian measure we need the following
\begin{lemma}  \label{l: Gaussian measure}
Let $x_1, \ldots x_s$ be vectors in $\R^s$. Then
\[
  \g_s \left (\sum_{j=1}^s [-x_j,x_j] \right )
  \le \g_s( M \cdot B_{\infty}^s),
\]
where $M= \max_{j=1, \ldots s} \|x_j\|_2$.
\end{lemma}

The sum in the Lemma is understood as the Minkowski sum of sets of vectors, 
$A+B = \{a+b \;|\; a \in A, \; b \in B\}$.  

\medskip

\proof Let $F= \Span (x_1, \ldots x_{s-1})$ and let $V=F^{\perp}$.
Let $v \in V$ be a unit vector. Set $Z= \sum_{j=1}^{s-1}
[-x_j,x_j]$. Then
\begin{align*}
   \g_s \Bigl(\sum_{j=1}^s [-x_j,x_j] \Bigr)
   &= \int_V \g_F \Bigl( \Bigl( \sum_{j=1}^s [-x_j,x_j]-tv \Bigr) \cap F
                  \Bigr) \, d \g_V(t) \\
   &= \int_{[-P_V x_s, P_V x_s]} \g_F (Z+ t P_F x_s) d \g_V(t).
\end{align*}
By Anderson's Lemma (see \cite{Lif}), 
$\g_F (Z+ t P_F x_s) \le \g_F (Z)$. Thus,
\[
  \g_s \Bigl( \sum_{j=1}^s [-x_j,x_j] \Bigr)
  \le \g_V([-P_V x_s, P_V x_s]) \cdot \g_F(Z)
  \le \g_1([-M,M]) \cdot \g_F(Z).
\]
The proof of the Lemma is completed by induction.
\endproof

The Gaussian measure of a projection of the cube can be estimated
as follows.
\begin{proposition}
  Let $H$ be any subspace in $G_{n,n-k}$, $k < n/2$.
  Then
  \begin{equation} \label{measure of proj of cube}
  \gamma_H \Bigl( \frac{c}{\sqrt{k}} \sqrt{\log \frac{n}{k}} \,
                       P_H B_\infty^n \Bigr)
  \le e^{-cn/k}.
  \end{equation}
\end{proposition}

\proof Decompose $I$ into the disjoint union of the sets $J_1,
\ldots J_{s+1}$, so that each of the sets $J_1, \ldots J_s$
contains $k+1$ elements and $(k+1)s<n \le (k+1)(s+1)$. Let $1 \le
j \le s$. Let $U_j = H \cap (P_He_i, \ i \in \{1, \ldots n\}
\setminus J_j)^{\perp}$, where $e_1, \ldots e_n$ is the standard
basis of $\R^n$. Then $U_j$ is a one-dimensional subspace of $H$.
Set
\[
  x_j= \sum_{i \in J_j} \e_i P_He_i,
\]
where the signs $\e_i \in \{-1,1\}$ are chosen to maximize
$\|P_{U_j}x_j\|_2$. Let $E= \Span (x_1, \ldots x_{s-1})$. Since
$P_{U_j} B_{\infty}^n = [-x_j,x_j]$, we get
\[
  P_H B_{\infty}^n \cap E = \sum_{j=1}^s [-x_j,x_j],
\]
where the sum is understood in the sense of Minkowski addition.
Since $\|P_{U_J}\| =1$, $\|x_j\|_2 \le C \sqrt{k}$ and by Lemma
\ref{l: Gaussian measure},
\[
  \gamma_E \left ( \frac{\bar{c}\sqrt{\log s}}{\sqrt{k}}
  \sum_{j=1}^s [-x_j,x_j] \right )
  \le \gamma_E ( c'\sqrt{\log s} \cdot B_{\infty}^E) \le e^{-cs}
\]
for some appropriately chosen constant $\bar{c}$. Finally,
log-concavity of the Gaussian measure implies that for any convex
symmetric body $K \subset H$
\[
  \gamma_H (K) \le \gamma_E(K \cap E).
\]
\endproof

Combining \eqref{upper P} and \eqref{measure of proj of cube} we
obtain $P \le 2e^{-cR}$, whenever $R \le c \log (m/r)$.

\subsection{Robustness and codes for finite alphabets}
Robustness is a well known property of the Basis Pursuit method.
It states that the solution to (BP) is stable with respect to the $1$-norm.
Indeed, it is not hard to show that, once Theorem \ref{ecc} holds,
the unknown vector $y$ in Theorem \ref{ecc} can be approximately recovered
from $y'' = y' + h$, where $h \in \R^m$ is any additional
error vector of small $1$-norm (see \cite{CT}).
Namely, the solution $u$ to the Basis Pursuit problem
$$
\min_{u \in Y} \|u - y''\|_1
$$
satisfies
$$
\|u - y\|_1 \le 4 \|h\|_1.
$$
This implies a possibility of quantization of the coefficients
in the process of encoding and yields {\em error correcting codes over
alphabets of size polynomial in $n$}.

The following is the $(m,n,r)$-error correcting code under
assumption \eqref{mnr'}, with input words $x$ over the alphabet
$\{1,\ldots,p\}$ and the encoded words $y$ over the alphabet
$\{1, \ldots, C p n^{3/2}\}$. The construction is the same as
in \eqref{ss:ecc}; we just introduce quantization.
The encoder takes $x \in \{1,\ldots,p\}^n$, computes
$y = Qx$ and outputs the $\hat{y}$ whose coefficients are the quantized
coefficients of $y$ with step $\frac{1}{10m}$.
Then $\hat{y} \in \frac{1}{10m} \Z^m \cap [-p\sqrt{m}, p\sqrt{m}]^m$,
which by rescaling can be identified with $\{1, \ldots, C p n^{3/2}\}$
because we can assume that $m \le 2n$.
The decoder takes $y' \in \frac{1}{10m} \Z^m$, finds solution $u$
to (BP) with $Y = \range(Q)$, inverts to $x' = Q^T u$ and
outputs $\hat{x'}$ whose coefficients are the quantized
coefficients of $x'$ with step $1$.

This is indeed an $(m,n,r)$-error correcting code. If
$y'$ differs from $\hat{y}$ on at most $r$ coordinates, this and
the condition $\|\hat{y} - y\|_1 \le \frac{1}{10}$ implies
by the robustness that $\|u-y\|_1 \le 0.4$. Hence
$\|x'-x\|_2 = \|Q^T (u-y)\|_2 = \|u-y\|_2 \le \|u-y\|_1 \le 0.4$.
Thus $\hat{x'} = x$, so the decoder recovers $x$ from $y'$ correctly.

The robustness also implies a ``continuity'' of our error correcting
codes. If the number of corrupted coordinates in the received message
$y'$ is bigger than $r$ but is still a small fraction,
then the $(m,n,r)$-error correcting code above can still recover $y$
up to some small fraction of the coordinates.

We hope to return to consequences of our method, in particular
to robustness and continuity of our codes and generally to codes over
finite alphabets, in a separate publication.

{\small

\end{document}